\catcode`\@=11

\font\tenmsb=msbm10
\font\sevenmsb=msbm7 \font\fivemsb=msbm5  \newfam\msbfam
\textfont\msbfam=\tenmsb\scriptfont\msbfam=\sevenmsb
\scriptscriptfont\msbfam=\fivemsb

\def\hexnumber@#1{\ifnum#1<10 \number#1\else \ifnum#1=10 A\else\ifnum#1=11
 B\else\ifnum#1=12 C\else \ifnum#1=13 D\else\ifnum#1=14 E\else\ifnum#1=15
 F\fi\fi\fi\fi\fi\fi\fi}
 \def\msb@{\hexnumber@\msbfam}

\mathchardef\hbar="0\msb@7E

\def\Bbb{\ifmmode\let\next\Bbb@\else
\def\next{\errmessage{Use \string\Bbb\space only in math mode}}\fi\next}
\def\Bbb@#1{{\Bbb@@{#1}}} \def\Bbb@@#1{\fam\msbfam#1}

\catcode`\@=\active

\def\del{\partial}
\def\CR{\hbox{{$\cal R$}}} \def\CM{\hbox{{$\cal M$}}}
\def\cg{\hbox{{\sl g}}} 
\def\cb{\hbox{{\sl b}}}
\def\cf{\hbox{{\sl f}}}

\def\C{{\Bbb C}}

\def\rbiprod{{\cdot\kern-.33em\triangleright\!\!\!<}}
\def\lbiprod{{>\!\!\!\triangleleft\kern-.33em\cdot\, }}
\def\tens{\mathop{\otimes}}
\def\la{{\triangleright}}
\def\extd{{{\rm d}}}

\def\Ad{{\rm Ad}}
\def\ad{{\rm ad}}
\def\ev{{\rm ev}}
\def\id{{\rm id}}
\def\<{\langle}
\def\>{\rangle}

\def\haj#1{{\mathaccent20 {#1}}}
\def\und#1{{\underline {#1}}}
\def\check{\haj{\ }}

\def\o{{}_{\scriptscriptstyle(1)}}
\def\t{{}_{\scriptscriptstyle(2)}}
\def\bo{{}^{\bar{\scriptscriptstyle(1)}}}
\def\bt{{}^{\bar{\scriptscriptstyle(2)}}}
\def\uo{{{}^{\scriptscriptstyle(1)}}}
\def\ut{{{}^{\scriptscriptstyle(2)}}}
\def\Bo{{{}_{\und{\scriptscriptstyle(1)}}}}
\def\Bt{{{}_{\und{\scriptscriptstyle(2)}}}}

\def\cmath#1{\[\begin{array}{c} #1 \end{array}\]}
\def\eqn#1#2{\begin{equation}#2\label{#1}\end{equation}}

\documentstyle[11pt]{article}
\begin{document}\baselineskip 22pt

{\ }\qquad  \hskip 4.3in DAMTP/96-83
\vspace{.2in}

\begin{center} {\LARGE BRAIDED GEOMETRY AND THE INDUCTIVE CONSTRUCTION OF LIE ALGEBRAS AND QUANTUM GROUPS}
\\ \baselineskip 13pt{\ }
{\ }\\
Shahn Majid\footnote{Royal Society University Research Fellow and Fellow of
Pembroke College, Cambridge}\\
{\ }\\
Department of Applied Mathematics \& Theoretical Physics\\
University of Cambridge, Cambridge CB3 9EW\\
+\\
Department of Mathematics, Harvard University\\
Science Center, Cambridge MA 02138, USA\footnote{During 1995+1996}\\
\end{center}
 
\vspace{10pt}
\begin{quote}\baselineskip 13pt
\noindent{\bf Abstract}
Double-bosonisation associates to a braided group in the
category of modules of a quantum group, a new quantum group. We announce the
semiclassical version of this inductive construction.

\end{quote}
\baselineskip 22pt

\section{Introduction}

A question usually overlooked in deformation theory is that of uniformity: we
can quantise this or that Poisson manifold, but do our individual quantisations
fit together into a coherent `quantum geometry'? In classical geometry the
co-ordinate rings are assumed {\em uniformly} to be commutative. When we relax
this, each object has many `directions' in which to become non-commutative and
we need to know how to pick these in a coherent way.

This problem is addressed  by braided geometry, introduced by the author
through about 60 papers since 1989. Rather than deforming one algebra at a
time, we deform the tensor product itself; we do group theory and geometry in a
braided category in place of ${\rm Vec}$. Then all mathematical concepts
founded in linear algebra are $q$-deformed uniformly as we switch on the
braiding. Braided geometry has its own method of proofs in which algebraic
information `flows' along braid and tangle diagrams like information in a
computer, except that under and over crossings of wires are nontrivial braiding
operators $\Psi$
\cite{M1}\cite{M2}\cite{M3}\cite{M3}\cite{M4}\cite{M6}\cite{M7}.
In physical
terms, braided geometry is a generalisation of supergeometry with $-1$ in
Bose-Fermi statistics replaced by braid statistics (e.g. by   $q$). This
is conceptually quite different from the usual quantisation picture where
$q=e^{\hbar\over 2}$. But braided-commutative with respect to some $\tens_q$
still means non-commutative with respect to the usual $\tens$, so we generate
noncommutative algebras, which we can then `semiclassicalise' via such an
expansion; we do not start with Poisson brackets but rather we generate them,
i.e. this is a deeper point of view. \goodbreak

The starting point is the concept of {\em braided group}\cite{M1} or Hopf
algebra in a braided category. This means an algebra and coalgebra $B$ in the
category for which the coproduct $\und\Delta:B\to B\und\tens B$ is an algebra
homomorphism, where $\und\tens$ is the {\em braided tensor product} of
algebras\cite{M1} in a braided category. In concrete terms, $B\und\tens B$
has product $(a\tens b)(c\tens d)=a\Psi(b\tens c)d$.

The simplest example\cite{M4} is the tensor algebra $TV$ on a
finite-dimensional vector space $V$ equipped with a braiding $\Psi:V\tens V\to
V\tens V$. Write $TV=\C\<x_i\>$ and $\Psi(x_i\tens x_j)=x_b\tens
x_aR^a{}_i{}^b{}_j$ (with summation of indices), where $R$ obeys the
Yang-Baxter equation. Then the coproduct has the form\cite{M4}:
\[ \und\Delta x_{i_1}x_{i_2}\cdots x_{i_m}=\sum_{r=0}^mx_{j_1}\cdots
x_{j_r}\tens x_{j_{r+1}}\cdots x_{j_m}\left[{m\atop r};R\right]^{j_1\cdots
j_m}_{i_1\cdots i_m}.\]
In his talk, Rosso\cite{Ros:bat} mentioned the `quantum shuffle algebra' but
this is just the graded dual of $TV$. Its product has just the structure of the
coproduct of the latter. Writing $y^{i_m\cdots i_1}$ for the dual basis to
$x_{i_1}\cdots x_{i_m}$, clearly
\[ y^{i_m\cdots i_{r+1}}\cdot y^{i_r\cdots i_1}=\left[{m\atop
r};R\right]^{i_1\cdots i_m}_{j_1\cdots j_m}\!\!\!y^{j_m\cdots j_1},\quad
\und\Delta y^{i_m\cdots i_1}=\sum_{r=0}^my^{i_m\cdots i_{r+1}}\tens
y^{i_r\cdots i_1}.\]
{}From standard properties\cite{M4} of these {\em braided binomial}
matrices $[{m\atop r};R]$,
\[\pi:T(V^*)\to (TV)^*,\quad \pi(y^{i_m}y^{i_{m-1}}\cdots
y^{i_i})=[m;R]!{}^{i_1\cdots i_m}_{j_1\dots j_m}y^{j_m\cdots j_1}\]
is a homomorphism of braided groups, where $[m;R]!$ are the {\em braided
factorial} matrices and $T(V^*)=\C\<y^i\>$. Hence $\ev:T(V^*)\tens TV\to \C$,
\eqn{ev}{\ev(f(y),g(x))=\pi(f(y))(g(x))=f(\del)g(x)|_{x=0}
=f(y)g(\overleftarrow{\del})|_{y=0}}
is a duality pairing of braided groups. Here $\del$ denotes {\em braided
differentiation}
\[\del^i x_{i_1}\cdots x_{i_m}=x_{j_2}\cdots x_{j_m}[m;R]^{ij_2\cdots
j_m}_{i_1\cdots i_m}\]
where $[m;R]$ is the {\em braided integer} matrix. Similarly for
$\overleftarrow{\del}$. These are some rudiments of braided geometry on free
algebras\cite{M4}.

The kernels of $\ev$ may be non-zero; quotienting by them gives new braided
groups such as the quantum planes $\C_q^n$. Another choice\cite{Ma:csta} is
$R^i{}_j{}^k{}_l=\delta^i{}_j\delta^k{}_lq^{\beta_{jl}}$, where $\beta$ is a
bilinear form. This is the case considered in \cite{Ros:bat} and Fronsdal's
talk\cite{Fro:gen}.
When $\beta$ comes from a Cartan matrix, Lusztig\cite{Lus} computed $\ker\pi$
as the $q$-Serre relations, i.e. $T(V^*)/\ker\pi={\rm image}(\pi)=U_q(n_+)$.

\section{Transmutation and Bosonisation; Induction Principle}
General theorems about braided groups are the following. We use Sweedler
notation $\Delta h=h\o\tens h\t$ for coproducts, $S$ for the antipode and
$\CR=\CR\uo\tens\CR\ut$ for Drinfeld's quasitriangular structure (summations
understood). \goodbreak

{\bf 1. Transmutation} (SM 1990). Let $H$ be a quantum group (quasitriangular
Hopf algebra). Its {\em transmutation} is the braided group $\und H\in{}_H\CM$,
the braided category of modules. $\und H$ is $H$ as a module-algebra by $\Ad$,
and
\[ \und\Delta h=h\o S\CR\ut\tens \Ad_{\CR\uo}(h\t),\quad \Psi(h\tens
g)=\Ad_{\CR\ut}(g)\tens\Ad_{\CR\uo}(h).\]

{\bf 2. Bosonisation} (SM 1991). Let $B\in {}_H\CM$ be a braided group with
$\und\Delta b=b\Bo\tens b\Bt$ and action $\la$ of $H$. Its  {\em bosonisation}
is the Hopf algebra $B\lbiprod H$ generated by $H$ as a Hopf algebra, $B$ as an
algebra, and
\eqn{bos}{ hb=(h\o\la b)h\t,\quad \Delta b=b\Bo\CR\ut\tens \CR\uo\la b\Bt.}

{\bf 3. Biproducts} (cf. Radford 1985, SM 1992). Let $B\in{}^{F}_{F}\CM$, the
crossed modules over a Hopf algebra $F$ with bijective $S$. There is a {\em
biproduct} Hopf algebra $B\lbiprod F$ projecting to $F$.
Every projection to $F$ is of this form.

{\bf 4. Double-bosonisation} (SM 1995). Let $B\check$ be dually paired to
$B\in{}_H\CM$ via $\ev:B\tens B\check\to \C$. There is a quantum group
$B\lbiprod H\rbiprod B\check^{\rm op}$ containing $B\lbiprod H$ and $H\rbiprod
B\check^{\rm op}$ as subHopf algebras, defined by (\ref{bos}) and
\cmath{ b\Bo \CR\ut c\Bo \ev(\CR\uo\la b\Bt, c\Bt)=\ev(b\Bo,\CR\ut\la
c\Bo)c\Bt\CR\uo b\Bt\\
hc=(h\t\la c)h\o,\quad \Delta c=\CR\ut\la c\Bo\tens c\Bt\CR\uo,\quad \CR^{\rm
new}=\CR\exp^{-1},}
where $\CR^{\rm new}$ needs a canonical element (coevaluation) $\exp\in
B\check\tens B$ for $\ev$.

{\bf 5. Double-biproducts} (SM 1995). Let $B\check\in{}^{F}_{F}\CM$ be dually
paired to $B\in{}^{F}_{F}\CM$ in 3. as in\cite{M6}. There is a Hopf algebra
$B\lbiprod
F\rbiprod B\check^{\rm op}$. A functor ${}_H\CM\to{}^H_H\CM$ allows
2. \& 4. to be viewed as special cases of 3. \& 5.

Bosonisation has been used to construct inhomogeneous Hopf algebras
$\C^n_q\lbiprod \widetilde{U_q(su_n)}$. The $\widetilde{\ }$ denotes a central
extension.
On the other hand, double-bosonisation can be iterated to provide a graph of
quantum groups, including the standard families of $U_q(\cg)$ as well as new
quantum groups without classical limit. At each node $H$, the branches are the
inequivalent $B\in {}_H\CM$. The new node is $B\lbiprod H\rbiprod B\check^{\rm
op}$. The initial node is the quantum group $\C$. Its central extension is the
quantum line $U_q(1)$. Adjoining the braided line $\C_q$ to this yields
$U_q(su_2)$. There are several braided groups in the category of
$\widetilde{U_q(su_2)}$-modules, each yielding a new quantum group. The
quantum-braided plane $\C_q^2$ gives us $U_q(su_3)$. There are some
technicalities, see \cite{M6}.

The required quantum-braided planes for induction up the A,B,C,D series
$U_q(\cg)$ are known, while the exceptional series are
currently under investigation. As there is surely {\em some}
braided group $B$ in the category of $\widetilde{U_q(e_8)}$-modules, we obtain
at least one quantum group $B\lbiprod \widetilde{U_q(e_8)}\rbiprod B\check^{\rm
op}$ which could be called $U_q(e_9)$! Presumably it does not survive as $q\to
1$. Also, building up $U_q(\cg)$ inductively by a series of triple products
yields automatically a natural {\em inductive block basis} for it, which
becomes a basis when we fix bases for the braided planes $B$ which are adjoined
at each stage.  For example,
\eqn{Uqsl}{
U_q(su_n)=\C_q^{n-1}\lbiprod\C_q^{n-2}\lbiprod\cdots\lbiprod\C_q\lbiprod
U_q(\beta)\rbiprod\C_q\rbiprod \cdots\rbiprod \C_q^{n-2}\rbiprod\C_q^{n-1}}
where the central extensions are collected together as $U_q(\beta)=U(1)^{\tens
n}$ generated by $H_i$ with a quasitriangular structure
$\CR_\beta=q^{\sum\beta^{-1}_{ij}H_i\tens H_j}$. Here $\beta$ is the
symmetrised Cartan matrix. This {\em proves} the PBW theorem for $U_q(\cg)$ and
explicitly constructs $U_q(n_+)=\C_q^{n-1}\lbiprod \C_q^{n-2}\cdots
\lbiprod\C_q$. Choosing bases for the $\C_q^i$ gives us a basis for $U_q(n_+)$,
as well as all the relations between them (including the $q$-Serre relations
when expressed in terms of the simple roots).
The inductive basis is coherent across the graph of quantum groups. Moreover,
its restriction to any substring of factors gives a sub-braided or quantum
group. In (\ref{Uqsl}), $\C_q^{n-1}\lbiprod\cdots\lbiprod \C_q\lbiprod
U_q(\beta)=U_q(b_+)$, $\C_q^2\lbiprod \C_q\lbiprod U_q(\beta)\rbiprod \C_q=
\C_q^2\lbiprod\widetilde{U_q(su_2)}$, etc. If one is interested in only half
the story, i.e. only in  constructing $U_q(b_+)$, one can also do it by
iterated biproducts. Thus, $\C_q^{n}\lbiprod\widetilde{U_q(b_+)}$ gives the
$q$-Borel of $U_q(su_{n+1})$.

Double-bosonisation also generalises Lusztig's construction. Any $\beta$
defines a quantum group $U_q(\beta)$ with generators $h_i$ and
$\CR_\beta=q^{\sum\beta_{ij}h_i\tens h_j}$. $B=\C\<y^i\>$, $B\check=\C\<x_i\>$,
paired by (\ref{ev}), live in the category of $U_q(\beta)$-modules by $h_i\la
y^j=\delta_{ij}y^j$, $h_i\la x_j=-\delta_{ij}x_j$. So $\C\<y^i\>\lbiprod
U_q(\beta)\rbiprod\C\<x_i\>$ is a Hopf algebra. Quotienting by the kernels of
$\ev$ we obtain a quantum group $U_q(n_+)\lbiprod U_q(\beta)\rbiprod U_q(n_-)$
with $\CR=\CR_\beta\exp^{-1}$.
For generic $\beta$ (or generic $R$-matrix in Section~1)
$[m;R]$ are invertible and the coevaluation for (\ref{ev}) is
\cmath{ \exp=\sum_{m=0}^{\infty}x_{i_m}\cdots x_{i_1}([m;R]!^{-1})^{i_1\cdots
i_m}_{j_1\cdots j_m}y^{j_m}\cdots y^{j_1}\in B\check^{\rm op}\tens B.}
Otherwise, quotienting by the kernels is nontrivial but we still have $\del$,
$\overleftarrow{\del}$ and the {\em braided exponential} $\exp$ is
characterised as their eigenfunction\cite{M7}.

Note that Fronsdal in his talk and \cite{Fro:gen} considered recursion
relations for an ansatz of the form $\CR_\beta f(x,y)$ to obey the Yang-Baxter
equation, with resulting Hopf algebra being coboundary. By contrast,
double-bosonisation already provides a closed expression for $\CR$ via a
braided-exponential, proves that quotienting by kernels of $\ev$ yields a Hopf
algebra and proves that it is quasitriangular. \cite{M6} has been
circulated in October 1995.

\section{Braided-Lie Bialgebras and Lie Induction}

We now announce a semiclassical concept of braided groups. Let $\cg$,
$\delta:\cg\to \cg\tens \cg$, $r\in \cg\tens \cg$ be a quasitriangular Lie
bialgebra as per Drinfeld\cite{Dri}. Let $2r_+=r+\tau(r)$ where $\tau$ is
transposition. Let $\la$ denote an action of $\cg$.\goodbreak

{\bf 0.} A {\em braided-Lie bialgebra} $\cb\in{}_{\cg}\CM$ is a $\cg$-covariant
Lie algebra and $\cg$-covariant Lie coalgebra with cobracket
$\und\delta:\cb\to \cb\tens \cb$ obeying $\forall x,y\in\cb$,
\[\und\delta([x,y])=\ad_x\delta y-\ad_y\delta x-\psi(x\tens y);\quad
\psi=2r_+(\la\tens\la)\circ (\id-\tau),\]
i.e., $\extd\und\delta=\psi$ where $\extd$ is the Lie coboundary on
$\und\delta\in C^1_{\ad}(\cb,\cb\tens\cb)$ and $\extd\psi\equiv0$.

{\bf 1.} Let $i:\cg\to\cf$ be a map of Lie bialgebras. The {\em transmutation}
of $\cf$ is a braided-Lie bialgebra $\und\cf\in{}_{\cg}\CM$  with Lie algebra
$\cf$ and for all $x\in \cf$,
\[ \und\delta x=\delta x + r\uo\la x\tens i(r\ut)-i(r\ut)\tens r\uo\la x,\quad
\la=\ad\circ i.\]
In particular, $\cg$ has a braided version $\und\cg\in{}_{\cg}\CM$ by $\ad$,
the same bracket, and
\eqn{kk}{\und\delta x= 2r_+\uo\tens [x,r_+\ut].}

{\bf 2.} Let $\cb\in{}_{\cg}\CM$ be a braided-Lie bialgebra. Its {\em
bosonisation} is the Lie bialgebra $\cb\lbiprod \cg$ with $\cg$ as sub-Lie
bialgebra, $\cb$ as sub-Lie algebra and
\eqn{liebos}{[\xi,x]=\xi\la x,\quad \delta x=\und\delta x+r\ut\tens r\uo\la
x-r\uo\la x\tens r\ut,\quad\forall \xi\in\cg,\ x\in\cb.}

{\bf 3.} Let $\cf$ be a Lie bialgebra and ${}_{\cf}^{\cf}\CM$ its category of
Lie crossed modules (=modules of the Drinfeld double $D(\cf)$.) Objects $\cb$
are simultaneously $\cf$-modules $\la$ and $\cf$-comodules
$\beta:\cb\to\cf\tens\cb$ obeying $\forall \xi\in\cf, x\in\cb$,
\[ \beta(\xi\la x)=([\xi, ]\tens\id+\id\tens \xi\la)\beta(x)+(\delta \xi)\la
x.\]
Writing $\beta(x)=x\bo\tens x\bt$, the infinitesimal braiding in this category
is
\[\psi(x\tens y)= y\bo\la x\tens y\bt-x\bo\la y\tens x\bt- y\bt\tens y\bo\la x+
x\bt\tens x\bo\la y.\]
Let $\cb\in{}_{\cf}^{\cf}\CM$ be a braided-Lie bialgebra. The {\em bisum} Lie
bialgebra $\cb\lbiprod \cg$ has semidirect Lie bracket/cobracket and projects
onto $\cf$. Any Lie bialgebra projecting onto $\cf$ is of this form. A functor
${}_{\cg}\CM\to{}^{\cg}_{\cg}\CM$ relates 2. \& 3.

{\bf 4.} Let $\cb\check\in{}_{\cg}\CM$ be a braided-Lie bialgebra dually paired
with $\cb$ by invariant $\ev:\cb\tens\cb\check\to \C$. Its {\em
double-bosonisation} is the Lie bialgebra $\cb\lbiprod\cg\rbiprod\cb\check^{\rm
op}$ with $\cg$ as sub-Lie bialgebra, $\cb,\cb\check^{\rm op}$ sub-Lie
algebras, (\ref{liebos}) and
\cmath{{}[\xi,\phi]=\xi\la \phi,\quad
[x,\phi]=\ev(x\Bo,\phi)x\Bt+\ev(x,\phi\Bo)\phi\Bt+2r_+\uo\ev(x,r_+\ut\la
\phi)\\
\delta\phi=\und\delta\phi+r\ut\la \phi\tens r\uo-r\uo\tens r\ut\la \phi,\quad
r^{\rm new}=r-\sum_a f^a\tens e_a,}
$\forall x\in\cb,\xi\in\cg$ and $\phi\in\cb\check$. Here $\und\delta
x=x\Bo\tens x\Bt$, etc., and $r^{\rm new}$ assumes that $\ev$ has a
coevaluation, i.e. if $\{e_a\}$ is a basis of $\cb$ then $\{f^a\}$ is dual
w.r.t. $\ev$.

Double bosonisation provides an inductive  construction for quasitriangular Lie
bialgebras, preserving factorisability (nondegeneracy of $r_+$). It is a
co-ordinate free version of the idea of adjoining a node to a Dynkin diagram
(adjoining a simple root vector in the Cartan-Weyl basis). Moreover, building
up $\cg$ iteratively like this also builds up the quasitriangular structure
$r$. The braided-Lie bialgebra used in the induction could be trivial:

{\bf Proposition.} Let $\cg$ be a semisimple factorisable (s.s.f) Lie bialgebra
 and $\cb$ a faithful isotypical representation such that $\Lambda^2\cb$ is
isotypical. Then $\cb$ with zero bracket and zero cobracket is a braided-Lie
bialgebra in ${}_{\tilde{\cg}}\CM$, and
$\cb\lbiprod\tilde{\cg}\rbiprod\cb\check^{\rm op}$ is another s.s.f. Lie
bialgebra. Here $\tilde{\cg}$ is a central extension.

The induction also works at the simple strictly quasitriangular level (with $b$
irreducible). For example, the $2$-dimensional and $3$-dimensional
representations of $su_2$ have the required property (ensuring $\psi\propto
(\id-\tau)$ in ${}_{\cg}\CM$). They give $su_3$ and $so_5$, taking us up the
$A$ and $B$ series respectively.

Finally, just as Lie bialgebras extend to Poisson-Lie groups, so braided-Lie
bialgebra structures generally extend to the associated Lie group of $\cg$. The
resulting Poisson bracket does not, however, respect the group product in the
usual way but rather up to a `braiding' obtained from $\psi$.

{\bf Example.} The transmutation (\ref{kk}) of the
Drinfeld-Sklyanin (or other factorisable) Lie cobracket on semisimple $\cg$ is
the
Kirillov-Kostant Lie cobracket. Moreover, this {\em Kirillov-Kostant
braided-Lie bialgebra} extends, in principle, to a braided-Poisson Lie group.
Details are to appear elsewhere.

\end{document}